\documentclass[reqno]{amsart}
\usepackage{amsfonts,graphicx,float}
\usepackage[arrow,matrix]{xy}

\restylefloat{figure}

\makeatletter                                                  %
\def\th@plain{\slshape}                                        %
\makeatother                                                   %

\newcommand{\tnorm}{\star}
\newcommand{\oi}{[0,1]}
\newcommand{\ooii}{\{0,1\}}

\newcommand{\Nbb}{\mathbb{N}}
\newcommand{\Zbb}{\mathbb{Z}}
\newcommand{\Qbb}{\mathbb{Q}}
\newcommand{\Rbb}{\mathbb{R}}
\newcommand{\two}{\mathbf{2}}
\newcommand{\Boole}{{\mathit{Boole}}}
\newcommand{\MV}{{\mathit{MV}}}

\newcommand{\llgroup}{$\ell$-group}
\newcommand{\llgroups}{$\ell$-groups}

\newcommand{\llhomomorphism}{$\ell$-homomorphism}
\newcommand{\llhomomorphisms}{$\ell$-homomorphisms}
\newcommand{\llisomorphism}{$\ell$-isomorphism}
\newcommand{\llisomorphic}{$\ell$-isomorphic}

\newcommand{\ffrak}{\mathfrak{f}}
\newcommand{\gfrak}{\mathfrak{g}}
\newcommand{\pfrak}{\mathfrak{p}}
\newcommand{\qfrak}{\mathfrak{q}}
\newcommand{\mfrak}{\mathfrak{m}}
\newcommand{\nfrak}{\mathfrak{n}}
\newcommand{\Nfrak}{\mathfrak{N}}
\newcommand{\Gfrak}{\mathfrak{G}}

\newcommand{\Ocal}{\mathcal{O}}
\newcommand{\Luk}{\L ukasiewicz}

\newcommand{\To}{\Rightarrow}

\newcommand{\newword}[1]{\textsl{#1}}
\newcommand{\operator}[1]{\mathbf{#1}}
\newcommand{\vect}[3]{#1_#2,\ldots ,#1_#3}
\newcommand{\projvect}[3]{#1_#2:\ldots :#1_#3}

\newcommand{\angles}[1]{\langle #1 \rangle}
\newcommand{\free}[2]{\Free_{#1}(#2)}

\newcommand{\class}[1]{\mathcal{#1}}

\DeclareMathSymbol{\upharpoonright}{\mathrel}{AMSa}{"16}
\let\restriction\upharpoonright
\DeclareMathSymbol{\nmid}{\mathrel}{AMSb}{"2D}

\DeclareMathOperator{\Spec}{Spec}

\DeclareMathOperator{\den}{den}
\DeclareMathOperator{\lHom}{\ell-Hom}
\DeclareMathOperator{\Free}{Free}

\DeclareMathOperator{\Rat}{Rat}
\DeclareMathOperator{\Fl}{F\ell}
\DeclareMathOperator{\FVL}{FVL}

\theoremstyle{plain}
\newtheorem{theorem}{Theorem}[section]
\newtheorem{lemma}[theorem]{Lemma}

\newtheorem{corollary}[theorem]{Corollary}

\theoremstyle{definition}
\newtheorem{definition}[theorem]{Definition}
\newtheorem{example}[theorem]{Example}

\newtheorem{addendum}[theorem]{Addendum}

\begin{document}

\bibliographystyle{plain}

\sloppy

\title[Logical substitutions]{Dynamical properties\\
of logical substitutions}

\author[G. Panti]{Giovanni Panti}
\address{Department of Mathematics\\
University of Udine\\
via delle Scienze 208\\
33100 Udine, Italy}
\email{panti@dimi.uniud.it}

\begin{abstract}
Many kinds of algebraic structures have associated dual topological spaces, among others commutative rings with $1$ (this being the paradigmatic example), various kinds of lattices, boolean algebras, $C^*$-algebras, \ldots. These associations are functorial, and hence algebraic endomorphisms of the structures give rise to continuous selfmappings of the dual spaces, which can enjoy various dynamical properties; one then asks 
about the algebraic counterparts of these properties.
We address this question from the point of view of algebraic logic. The datum of a set of truth-values and a ``conjunction'' connective on them determines a propositional logic and an equational class of algebras. The algebras in the class have dual spaces, and the duals of endomorphisms of free algebras provide dynamical models for Frege deductions in the corresponding logic.
\end{abstract}

\keywords{algebraic logic, free algebras, spectral spaces, dual mappings}

\thanks{\emph{2000 Math.~Subj.~Class.}: 03B50; 37A05}

\maketitle

\section{Introduction}

Everybody knows the classical truth-tables
$$
\begin{array}{c|cc}
\land & 0 & 1 \\
\hline
0 & 0 & 0 \\
1 & 0 & 1
\end{array}
\hspace{1cm}
\begin{array}{c|cc}
\lor & 0 & 1 \\
\hline
0 & 0 & 1 \\
1 & 1 & 1
\end{array}
\hspace{1cm}
\begin{array}{c|cc}
\to & 0 & 1 \\
\hline
0 & 1 & 1 \\
1 & 0 & 1
\end{array}
\hspace{1cm}
\begin{array}{c|c}
&\neg \\
\hline
0 & 1 \\
1 & 0
\end{array}
$$
and uses them automatically. The \newword{classical propositional calculus} studies the set of formulas that, when evaluated according to the truth-tables, always assume value~$1$.
One proceeds as follows:
\begin{enumerate}
\item a \newword{formula} is a polynomial built up from the propositional variables $x_i$ using the connectives $\land,\lor,\to,\neg,0,1$;
\item a \newword{valuation} is a function $p$, distributing over the connectives, from the set of formulas to $\ooii$;
\item a formula $r$ is \newword{true} if $p(r)=1$ for every valuation $p$;
\item a formula $r$ is \newword{deducible} if:
\begin{itemize}
\item[(a)] either is an element of a certain set $\Theta$ of basic axioms,
\item[(b)] or there exists a formula $s$ such that $s$ and $s\to r$ are deducible,
\item[(c)] or there exists a deducible formula $s$, propositional variables $\vect x1n$, and formulas $\vect t1n$, such that $r$ results from $s$ by substituting every $x_i$ that occurs in $s$ with the corresponding $t_i$;
\end{itemize}
\item the completeness theorem holds: a formula is true iff it is deducible.
\end{enumerate}

The completeness theorem relates a semantical notion ``the statement $r$ holds, regardless of the state of affairs $p$'' with a computational notion ``the statement $r$ can be deduced from certain statements using certain rules''. There are several computational procedures for which the completeness theorem holds: the ones sketched in~(4) are known as \newword{substitutional Frege systems}, and are the strongest ---in terms of minimizing the number of steps required to prove a true statement--- available proof systems~\cite{cookrec79}.
The two rules (4b) of \newword{Modus Ponens} and (4c) of \newword{substitution} have different flavors. The first rule is, in some sense, statical: if something is known ``locally'', i.e., concerns certain propositional variables, then conclusions are drawn involving the same variables. On the other hand, the substitution rule adds dynamics to the picture: local knowledge can be moved around. This is of course just a vague heuristic, and in the course of these notes we will give a precise formal ground to it.

We will work at a level of generality broader than that of classical logic, enlarging the set of truth-values to include more than \emph{true} and \emph{false}; such logical systems are known as \newword{many-valued logics}. Many-valued logic is an old discipline, going back to the twenties, and has recently been relived as a founding basis for fuzzy logic and fuzzy control; see~\cite{hajek98}, \cite{CignoliOttavianoMundici00}, \cite{gottwald01} for detailed presentations and further references.

The key ideas of this work are the following: given a set of truth-values $M\supseteq\ooii$, we introduce on it an algebraic structure, determined by the choice of a truth-table for the conjunction connective. We then consider the class $\operator{V}M$ of all algebras that are generated by $M$ in the sense of Universal Algebra, and we functorially associate a dual topological space to each object in $\operator{V}M$. Algebraic endomorphisms of certain objects of $\operator{V}M$ (the so-called free algebras) correspond to applications of the substitution rule in deductions in the logic determined by $M$. Moreover, such endomorphisms give rise to continuous selfmappings of the dual topological spaces. Any set $\Theta'\supseteq\Theta$ ($\Theta$ is a set of basic axioms as in~(4a)) is associated to an open set $O_{\Theta'}$ in the dual, and the deduction of new formulas from $\Theta'$ corresponds to taking the union of the backwards translates of $O_{\Theta'}$ under the dynamics.
Dynamical properties such as minimality or mixing have then logical consequences (see, e.g., Theorem~\ref{ref10}, Theorem~\ref{ref14}, and the discussion following Theorem~\ref{ref23}). It is worth remarking that the trade between the logical and the dynamical side may be beneficial to both: as an example, we obtain in Theorem~\ref{ref22} an intrinsic characterization of the differential of a piecewise-linear mapping, a concept introduced in~\cite{Tsujii01}.

A rather delicate point in our approach is the determination of the level of generality one should allow. Here we must really strike a balance: the stronger is the system (i.e., the more restrictions we put on $M$), the stronger are the results we obtain, and the more limited is the scope of the theory. The extreme case is in taking $M=\{0,1\}$, in which everything boils down to the Stone Duality. On the other extreme, one might relax the assumptions on $M$ to a bare minimum, even allowing cases in which the values $0$ and $1$ do not have a distinguished status: the only essential requirement seems to be that $\operator{V}M$ is a congruence-modular equational class. Of course, working at this level of generality requires a greater technical apparatus, and yields not easily visualizable results.

We stroke our balance by forcing $M$ to be a subset of the real unit interval $\oi$, and by insisting that the conjunction connective meets some natural restrictions. In the few places where we might have wished more elbow-room, we have added some Addenda to provide references for further developments. These Addenda are meant for people having some knowledge of Universal Algebra and lattice-ordered abelian groups, and may be safely skipped by the other readers.

\section{Many-valued logic}\label{ref24}

A \newword{t-norm} is a continuous function $\tnorm$ from $\oi^2$ to $\oi$ such that $(\oi,\tnorm,1)$ is a commutative monoid for which $a\le b$ implies $c\tnorm a\le c\tnorm b$. We have $a\tnorm 0=0$ for every $a$, since $a\le 1$ implies $0\tnorm a\le 0\tnorm 1=0$.
Every t-norm induces a binary operation $\to$ on $\oi$ via
$$
a\to b = \sup\{c:c\tnorm a\le b\}.
$$
Since $\tnorm$ is continuous, the defining $\sup$ is really a $\max$. We call $\to$ the \newword{implication} (or the \newword{residuum}) induced by $\tnorm$. One checks easily that the usual lattice operations on $\oi$ are definable from $\tnorm$ and $\to$ via $a\land b=a\tnorm(a\to b)$ and $a\lor b=\bigl((a\to b)\to b\bigl)\land
\bigl((b\to a)\to a\bigl)$. We also define $\neg a=a\to 0$.

The idea underlying these definitions is that $\tnorm$ is a function on truth-values representing a ``conjunction'' operator. Once a conjunction has been fixed, it is natural to define the truth-value of the implication $a\to b$ as the weakest value~$c$ such that the truth of the conjunction of $a$ and $c$ forces the truth of $b$. Note that ``weakest'' means ``truest'', i.e., nearest to $1$: one should regard a more implausible assertion as a stronger one. The above interrelationship of $\tnorm$ and $\to$ is usually expressed by saying that they constitute an \newword{adjoint pair}.

\begin{example}\label{ref3}
\begin{enumerate}
\item $a\tnorm b=a\land b$. One computes that
$$
a\to b=
\begin{cases}
1, & \text{if $a\le b$;} \\
b, & \text{otherwise;}
\end{cases}
\quad
\neg a=
\begin{cases}
1, & \text{if $a=0$;} \\
0, & \text{otherwise.}
\end{cases}
$$
This t-norm is usually called the \newword{G\"odel-Dummett conjunction}.
\item $a\tnorm b=ab$ (i.e., the ordinary product of $a$ and $b$). This is the \newword{product conjunction}, and we have
$$
a\to b=
\begin{cases}
1, & \text{if $a\le b$;} \\
b/a, & \text{otherwise;}
\end{cases}
\quad
\neg a=
\begin{cases}
1, & \text{if $a=0$;} \\
0, & \text{otherwise.}
\end{cases}
$$
\item $a\tnorm b=\max(a+b-1,0)$. Then
$$
a\to b=
\begin{cases}
1, & \text{if $a\le b$;} \\
1-(a-b), & \text{otherwise;}
\end{cases}
\quad
\neg a=1-a.
$$
These are the \newword{\Luk\ conjunction}, \newword{implication}, and \newword{negation}
\end{enumerate}
\end{example}

The above examples are in some sense exhaustive: by~\cite{MostertShields57} every t-norm is obtainable as a combination of these three basic t-norms.

Fix a cardinal number $\kappa$, either finite or countable, and 
define the set of propositional variables to be $\{x_i:i<\kappa\}$ (then either $\kappa=n$ and the propositional variables are $\vect x0{{n-1}}$, or $\kappa=\omega$ and the propositional variables are indexed by the natural numbers). Let $FORM_\kappa$ be the smallest set containing all propositional variables having index $<\kappa$, the constants $0$ and $1$, and such that, if $r,s\in FORM_\kappa$, then $(r\tnorm s),(r\to s)\in FORM_\kappa$.
A \newword{formula} is an element $r$ of $FORM_\omega=\bigcup_{n<\omega}FORM_n$. We sometimes write $r(x_{i_1},\ldots,x_{i_n})$ to signify that all propositional variables occurring in $r$ are among $x_{i_1},\ldots,x_{i_n}$. We drop parentheses according to the usual conventions, and we write $r\land s$, $r\lor s$, and $\neg r$ as abbreviations for 
$r\tnorm(r\to s)$, $\bigl((r\to s)\to s\bigl)\land
\bigl((s\to r)\to r\bigl)$, and $r\to 0$, respectively.

An \newword{algebra} is a set $A$ on which two binary operations $\tnorm_A,\to_A:A^2\to A$ and two elements $0_A,1_A\in A$ have been fixed.
Given a formula $r(\vect x0{{n-1}})$ and elements $\vect a0{{n-1}}\in A$, we write $r(\vect a0{{n-1}})$ for the element of $A$ obtained by replacing every $x_i$ with the corresponding $a_i$, and every operation symbol in $r$ with its realization in $A$ (the reader can easily supply a formal recursive definition). Given two formulas $r(\vect x0{{n-1}})$ and $s(\vect x0{{n-1}})$, we say that the identity $r=s$ is \newword{true} in $A$, and we write $A\models r=s$, if for every $\vect a0{{n-1}}\in A$ the elements $r(\vect a0{{n-1}})$ and $s(\vect a0{{n-1}})$ are equal.

\begin{example}\label{ref4}
\begin{enumerate}
\item Every singleton can be given the structure of an algebra in a unique trivial way; every identity is true in such an algebra.
\item Let $A=\oi$, endowed with the G\"odel-Dummett conjunction and implication, as in Example~\ref{ref3}(1). Then
$A\not\models\neg\neg x_0 = x_0$, so the 
double negation rule fails for the G\"odel-Dummett connectives (and analogously for the product connectives). On the other hand, $\neg\neg x_0 = x_0$ holds true in $\oi$ endowed with the \Luk\ connectives.
\end{enumerate}
\end{example}

Let $A,B$ be algebras. A mapping $\varphi:A\to B$ is a \newword{homomorphism} if it commutes with the connectives (i.e., $\varphi(a\tnorm_A b)=\varphi(a)\tnorm_B\varphi(b)$, $\varphi(0_A)=0_B$, and so on; in the following we will drop the subscripts). $A$ is [isomorphic to] a \newword{subalgebra} of $B$ if there exists an injective homomorphism from $A$ to $B$.
Let $\{A_j:j\in J\}$ be a family of algebras. The \newword{direct product} of the family is the algebra whose base set is the cartesian product $\prod_jA_j$, and in which the operations are defined componentwise; if all factors are equal, say to $A$, then we write $A^J$. If $\varphi:A\to B$ is a homomorphism, then the \newword{epimorphic image} $\varphi[A]$ of $A$ is a subalgebra of $B$.

Let $\class{A}$ be a class of algebras; then $\operator{H}\class{A}$ (respectively, $\operator{S}\class{A}$ and $\operator{P}\class{A}$) is the class of all epimorphic images (respectively, subalgebras and direct products) of algebras in~$\class{A}$. Note that we always work up to isomorphism, so we tacitly close every class we consider under isomorphic images.

\begin{definition}
A \newword{truth-value algebra} is a subalgebra $M$ of some algebra $A$ of the form $A=(\oi,\tnorm,\to,0,1)$, where $\tnorm$ and $\to$ are a t-norm and its residuum.
\end{definition}

Truth-value algebras are our basic building blocks.

\begin{example}
\begin{enumerate}
\item The set $\ooii$ is always closed under the operations, regardless of the specific t-norm we choose. Moreover, all t-norms induce the same structure on $\ooii$, namely that of the \newword{two-element boolean algebra}, which we denote by $\two$.
\item $M=\{0,1/m,2/m,\ldots,(m-1)/m,1\}$ endowed either with the \Luk\ connectives or the G\"odel-Dummett ones.
\end{enumerate}
\end{example}

For any class $\class{A}$ of algebras, let $\operator{V}\class{A}$ be the \newword{equational class} generated by $\class{A}$, i.e., the
class of all algebras in which are true all identities true in all algebras of $\class{A}$. More explicitly, the algebra $B$ is in $\operator{V}\class{A}$ iff, for every $r,s\in FORM_\omega$, if $A\models r=s$ for every $A\in\class{A}$, then $B\models r=s$. 
Garrett Birkhoff's completeness theorem~\cite[Theorem~II.11.9]{burrissan81} says that $\operator{V}\class{A}$ coincides with the class $\operator{HSP}\class{A}$ of all epimorphic images of subalgebras of products of algebras in~$\class{A}$.

We will consider classes of algebras of the form $\operator{V}M=\operator{HSP}M$, where $M$ is a truth-value algebra. 
We shall be concerned with two main cases:
\begin{itemize}
\item $\Boole=\operator{V}\two$. Elements of $\Boole$ are called \newword{boolean algebras};
\item if $M=\oi$ endowed with the \Luk\ connectives, then the elements of $\operator{V}M$ are called \newword{MV-algebras} (MV stands for \newword{Many-Valued}: the name is slightly misleading, since many-valued logic is not exhausted by \Luk\ logic, but it is firmly established; we accordingly write $\MV$ for the equational class
$\operator{V}M$).
\end{itemize}

A boolean algebra can be equivalently defined as a structure $A=(A,\land,\lor,\neg,0,1)$ such that
\begin{gather*}
x\land 1=x\lor 0=x;\\
x\land\neg x=0;\quad x\lor\neg x=1;\\
\text{$\land$ and $\lor$ are commutative and mutually distributive}.
\end{gather*}
Apart from the trivial change in the language ($\lor$ replaces $\to$), there is a theorem hidden in this equivalence, namely the fact that the above identities imply all other identities that hold in $\two$~\cite[p.~5]{Halmos63}.

An analogous alternative characterization of MV-algebras is obtained by adding a new connective $\oplus$ to the basic set 
$(\tnorm,\to,0,1)$.
We define $a\oplus b=\neg a\to b$, and directly compute that
$\oplus$ is \newword{truncated addition} on $\oi$, i.e., $a\oplus b=\min(a+b,1)$. Note that the basic set of connectives is equivalent to the set $(\oplus,\neg,0,1)$, since $a\tnorm b=\neg(\neg a\oplus\neg b)$ and $a\to b=\neg a\oplus b$. 
Then, in terms of the new set, an MV-algebra is a structure $(A,\oplus,\neg,0,1)$ 
such that $(A,\oplus,0)$ is an abelian monoid and the identities
$\neg\neg x=x$, $x\oplus 1=1$, $\neg(\neg x\oplus y)\oplus
y=\neg(\neg y\oplus x)\oplus x$ are satisfied~\cite[\S2]{mundicijfa}, \cite{CignoliOttavianoMundici00}.

\begin{lemma}
Let\label{ref5} $M$ be a truth-value algebra, $A\in\operator{V}M$. Then:
\begin{itemize}
\item[(i)] the operations $\land,\lor$ induce a lattice structure on $A$, with bottom element $0$ and top $1$;
\item[(ii)] the lattice order in~(i) is given by $a\le b$ iff $a\land b=a$ iff $a\to b=1$;
\item[(iii)] $A\models r=s$ iff $A\models (r\to s)\land(s\to r)=1$.
\end{itemize}
\end{lemma}
\begin{proof}
A structure $(A,\land,\lor,0,1)$ is a lattice with bottom and top iff it satisfies a certain finite set of identities (see, e.g.,~\cite[p.~28]{burrissan81}).
Since $M$ is totally-ordered, these identities are satisfied in $M$, and hence in $A\in\operator{V}M$. The first equivalence in~(ii) is just the definition of the lattice order on $A$. By definition of $\to$ in $M$, the identity $(x_0\land x_1)\to x_1=1$ is true in $M$, and hence in $A$. Therefore, if $a\land b=a$, then $a\to b=(a\land b)\to b=1$. On the other hand, if $a\to b=1$, then
$a\land b=a\tnorm(a\to b)=a\tnorm 1=a$.
This proves~(ii), and~(iii) is then immediate.
\end{proof}

By Lemma~\ref{ref5}(ii) we can deal with the ``less than'' relation
between formulas, thus writing
$A\models r(\vect x0{{n-1}})\le s(\vect x0{{n-1}})$
for $A\models r\to s=1$; this just means that however we choose $\vect a0{{n-1}}\in A$ we have $r(\vect a0{{n-1}})\le s(\vect a0{{n-1}})$. We then say that $r\le s$ is \newword{true} in $A$.

\begin{lemma}
Under\label{ref6} the same hypothesis as in Lemma~\ref{ref5}, the following relations are true in $A$:
\begin{itemize}
\item[(i)] $x_0\tnorm x_1\le x_0\land x_1$;
\item[(ii)] $x_0\le x_1\to(x_0\tnorm x_1)$;
\item[(iii)] $(x_0\to x_1)\tnorm(x_1\to x_2)\le x_0\to x_2$;
\item[(iv)] $(x_0\to x_1)\tnorm(x_2\to x_3)\le (x_0\tnorm x_2)\to
(x_1\tnorm x_3)$;
\item[(v)] $(x_0\to x_1)\tnorm(x_2\to x_3)\le (x_1\to x_2)\to
(x_0\to x_3)$.
\end{itemize}
\end{lemma}
\begin{proof}
One just checks that for every t-norm with residuum $\to$ the above relations are true in $(\oi,\tnorm,\to,0,1)$. Hence they are true in $M$, and therefore in every algebra in $\operator{V}M$.
\end{proof}

Fix now $\kappa$ and a truth-value algebra $M$. We want to construct an algebra $A$ in $\operator{V}M$ satisfying the following properties:
\begin{itemize}
\item $A$ is generated by a family $\{a_i:i<\kappa\}$ of elements indexed by $\kappa$;
\item if $r(x_{i_1},\ldots,x_{i_n})\in FORM_\kappa$ is not true in $M$, then the element $r(a_{i_1},\ldots,a_{i_n})\in A$ is different from $1$.
\end{itemize}
Essentially, this means that the $a_i$'s satisfy only those algebraic relations they cannot avoid, namely those that hold in $M$. Therefore they behave ``as freely as possible'', whence the name \newword{free algebra in $\operator{V}M$ over $\kappa$ generators} for $A$. 
Such an algebra is unique up to isomorphism, and can be characterized by an appropriate universal property: see,
e.g.,~\cite[II \S10]{burrissan81}. We write $\free \kappa{\operator{V}M}$ for $A$, and we construct it as follows: consider first $M^\kappa$, and let $a_i:M^\kappa\to M$ be the $i$-th projection. The $a_i$'s are elements of the algebra $M^{(M^\kappa)}$, and we define $\free \kappa{\operator{V}M}$ to be the subalgebra 
of $M^{(M^\kappa)}$ generated by them; the first condition is then automatically met. Suppose $M\not\models r$; then there exist
elements $b_{i_1},\ldots,b_{i_n}\in M$ such that $r(b_{i_1},\ldots,b_{i_n})\not=1$. Choose an element $c\in M^\kappa$ such that $a_i(c)=b_i$ for every $i\in\{\vect i1n\}$.
Then the projection of $r(a_{i_1},\ldots,a_{i_n})\in M^{(M^\kappa)}$ onto the $c$-th component has value 
$r(a_{i_1}(c),\ldots,a_{i_n}(c))=r(b_{i_1},\ldots,b_{i_n})
\not=1$: therefore $r(a_{i_1},\ldots,a_{i_n})$ is different from~$1$ in $M^{(M^\kappa)}$, and hence in $\free \kappa{\operator{V}M}$.

\begin{example}
Although the above construction looks baroque, it really works trivially. Suppose, e.g., we want to construct $\Free_3\Boole$. We first construct $\two^3$, which contains the eight elements $c_1=(0,0,0)$, $c_2=(0,0,1)$, \ldots, $c_8=(1,1,1)$. Then we construct $\two^{(\two^3)}$, which contains $2^8$ elements; three of these elements, namely
\begin{align*}
a_1 &= (0,0,0,0,1,1,1,1),\\
a_2 &= (0,0,1,1,0,0,1,1),\\
a_3 &= (0,1,0,1,0,1,0,1),
\end{align*}
correspond to the canonical projections to the first, second, and third component of the $c_j$'s (of course, the above explicit form for the $a_i$'s depends on how we listed the $c_j$'s). It is clear that if $r(x_1,x_2,x_3)$ is not $1$ in $\two$ for some choice of elements, then $r(a_1,a_2,a_3)\not=1=(1,1,1,1,1,1,1,1)$ in
$\free 3\Boole$: we just tried all possible choices!
\end{example}

Often it is not trivial to determine, for a given $M$ and $\kappa$, which are the elements of $\free\kappa{\operator{V}M}$, i.e., which functions from $M^\kappa$ to $M$ are expressible as polynomials over the projections $a_i$. The case that most concerns us is the \Luk\ one, which we will treat in Theorem~\ref{ref15}.
For the classical logic case, the answer is the following: give $\two$ the discrete topology, and $\two^\kappa$ the product topology. Then:
\begin{itemize}
\item[(i)] an element $f\in\two^{(\two^\kappa)}$ is in $\free \kappa\Boole$ iff it is continuous as a function $f:\two^\kappa\to\two$.
\end{itemize}
Since continuous functions from a topological space $X$ to $\two$ correspond to clopen subsets of $X$, this amounts to saying that the clopen subsets of $\two^\kappa$ are exactly the boolean combinations of the sets of the form $a_i^{-1}[1]=\{p\in\two^\kappa:p_i=1\}$. With this hint, we leave the proof of~(i) as an exercise for the reader. As a corollary we obtain:
\begin{itemize}
\item[(ii)] if $\kappa=n$ is finite, then $\two^n$ is a discrete space, and all functions $:\two^n\to\two$ are in $\free n\Boole$, i.e., are expressible by $n$-variable formulas. This is sometimes called the \newword{functional completeness} of the boolean connectives;
\item[(iii)] if $\kappa=\omega$ is countably infinite, then $\free\omega\Boole$ is the boolean algebra of all clopen subsets of the \newword{Cantor space} $\two^\omega$, the latter being the only compact, totally disconnected, second countable space having no isolated points~\cite[\S2.15]{hockingyou}.
\end{itemize}

\section{Spectral spaces}\label{ref13}

In the preceding section we have defined the equational class $\operator{V}M$ generated by a truth-value algebra $M$, and described the algebras $\free\kappa{\operator{V}M}$. In this section we functorially associate a dual topological space to each algebra in $\operator{V}M$; the duals of the free algebras are our main object of study.

We fix a truth-value algebra $M$; all algebras we consider are elements of $\operator{V}M$.
A \newword{filter} on $A\in\operator{V}M$ is the counterimage of $1$ under some homomorphism of domain $A$: sometimes filters are called \newword{dual ideals} since ideals, as in ring theory, are counterimages of $0$.

\begin{lemma}
A subset $\ffrak$ of $A$ is a filter iff it contains $1$ and is closed under Modus Ponens ($a,a\to b\in\ffrak$ implies $b\in\ffrak$). Every filter is closed under $\tnorm$ and $\land$, and is upwards closed ($b\ge a\in\ffrak$ implies $b\in\ffrak$). Given two homomorphisms $\varphi:A\to B$ and $\psi:A\to C$, we have $\varphi^{-1}[1]=\psi^{-1}[1]$ iff there exists an isomorphism $\chi:\varphi[A]\to\psi[A]$ such that $\psi=\chi\circ\varphi$.
\end{lemma}
\begin{proof}
Let $\ffrak=\varphi^{-1}[1]$; then clearly $1\in\ffrak$. If
$a,a\to b\in\ffrak$, then $\varphi(b)\ge\varphi(a\land b)=\varphi(a\tnorm(a\to b))=\varphi(a)\tnorm\varphi(a\to b)=1$. Conversely, assume that $\ffrak$ is a subset of $A$ containing $1$ and closed under Modus Ponens. Then $\ffrak$ is upwards closed, since $a\le b$ implies $a\to b=1\in\ffrak$, and hence $a\in\ffrak$ implies $b\in\ffrak$. If $a\in\ffrak$, then $b\to(a\tnorm b)\in\ffrak$ (by Lemma~\ref{ref6}(ii)); hence, if $b\in\ffrak$ as well, then $a\tnorm b,a\land b\in\ffrak$ (by Lemma~\ref{ref6}(i)).
Let now $\ffrak$ be a subset of $A$ containing $1$ and closed under MP.
Define a relation $\sim$ on $A$ by $a\sim b$ iff $a\to b,b\to a\in\ffrak$. Then $\sim$ is an equivalence relation (transitivity follows from Lemma~\ref{ref6}(iii)) which respects the operations. Indeed, if $a\sim b$ and $c\sim d$, then $a\tnorm c\sim b\tnorm d$ (by Lemma~\ref{ref6}(iv)) and $a\to c\sim b\to d$ (by Lemma~\ref{ref6}(v)). We can then form the quotient algebra $A/\ffrak$ in the natural way. If $a/\ffrak$ denotes the equivalence class of $a$ w.r.t.~$\sim$, then the map $\rho(a)=a/\ffrak$ is a surjective homomorphism from $A$ to $A/\ffrak$. It is then straightforward
to check that the map $\tau:A/\ffrak\to\varphi[A]$ defined by $\tau(a/\ffrak)=\varphi(a)$ is an isomorphism. Since $\varphi=\tau\circ\rho$, our last claim follows by composing isomorphisms.
\end{proof}

A filter $\pfrak$ is \newword{prime} if it is proper (i.e., different from $A$) and for every two filters $\ffrak,\gfrak$, if $\pfrak=\ffrak\cap\gfrak$ then either $\pfrak=\ffrak$ or $\pfrak=\gfrak$. A filter is \newword{maximal} if it is proper and not properly contained in any proper filter; clearly every maximal filter is prime. Although not difficult, the proof of the following lemma requires some knowledge of Universal Algebra; the reader can find a proof in~\cite[Proposition~1.3]{pantigeneric}.

\begin{lemma}
The\label{ref7} following are equivalent:
\begin{enumerate}
\item[(1)] $\pfrak$ is prime;
\item[(2)] $A/\pfrak$ is totally-ordered;
\item[(3)] $\pfrak=\varphi^{-1}[1]$, for some homomorphism $\varphi$ from $A$ to a totally-ordered algebra;
\item[(4)] the set of filters $\supseteq\pfrak$ is totally-ordered by inclusion;
\item[(5)] every filter $\supseteq\pfrak$ is prime;
\item[(6)] if $a\lor b\in\pfrak$, then either $a\in\pfrak$ or $b\in\pfrak$.
\end{enumerate}
\end{lemma}

Note that the only totally-ordered boolean algebra is $\two$ (if $a$ belongs to the totally-ordered boolean algebra $A$, then either $a\le\neg a$ or $\neg a\le a$, hence either $a\to\neg a=1$ or $\neg a\to a=1$; in the first case $a=0$, and in the second $a=1$), and therefore prime filters coincide with maximal ones. This simple fact distinguishes in a crucial way boolean algebras from MV-algebras and other algebras related to many-valued logics, as we will see later.

\begin{definition}
Let $A$ be an algebra, and let $\Spec A$ be the set of all prime filters of $A$. For every $a\in A$, let $O_a$ be the set of all $\pfrak\in\Spec A$ such that $a$ does not belong to $\pfrak$. Impose on $\Spec A$ the weakest topology in which all $O_a$'s are open (i.e., take the family of all $O_a$'s as an open subbasis). This is called the hull-kernel topology on $\Spec A$, and the resulting space is the \newword{spectral} (or \newword{dual}) \newword{space} of $A$.
\end{definition}

For every subset $D$ of $A$, let $O_D=\bigcup\{O_a:a\in D\}=
\{\pfrak:D\not\subseteq\pfrak\}$: it is an open set, and we will see in Theorem~\ref{ref2}(ii) that every open set has this form. We write $F_a=(\Spec A)\setminus O_a$ and $F_D=(\Spec A)\setminus O_D$
for the corresponding closed sets.

The mapping $A\mapsto\Spec A$ is functorial. Indeed, let $\varphi:A\to B$ be any homomorphism, and define $\varphi^*:\Spec B\to\Spec A$ by $\varphi^*(\pfrak)=\varphi^{-1}[\pfrak]$. 
$\varphi^*(\pfrak)$ is a prime filter because, if $\pfrak=\psi^{-1}[1]$ for some homomorphism $\psi$ from $B$ to a totally-ordered algebra
$C$, then $\varphi^*(\pfrak)$ is the kernel of $\psi\circ\varphi$, and the epimorphic image $(\psi\circ\varphi)[A]$ is totally-ordered, since it is a subalgebra of $C$.
We have 
$(\varphi^*)^{-1}[O_D]=\{\pfrak\in\Spec B:\varphi^*(\pfrak)\in O_D\}=\{\pfrak:D\not\subseteq\varphi^{-1}[\pfrak]\}=\{\pfrak:\varphi[D]\not\subseteq\pfrak\}=O_{\varphi[D]}$, and hence
$\varphi^*$ is continuous.
Moreover, $(\psi\circ\varphi)^*=\varphi^*\circ\psi^*$, so $\Spec$ is a contravariant functor from $\operator{V}M$ (viewed as a category with the homomorphisms as arrows) to the category of topological spaces and continuous mappings. 

\begin{theorem}
Let\label{ref2} $A$ be an algebra.
\begin{itemize}
\item[(i)]
The open sets in $\Spec A$ are in 1--1 correspondence with the filters of $A$, and this correspondence is an isomorphism w.r.t.~the $\subseteq$ relation.
\item[(ii)] We have $O_a\cap O_b=O_{a\lor b}$ and $O_a\cup O_b=O_{a\land b}$. The defining subbasis is intersection-closed, and an open set  is compact iff it is of the form $O_a$. $\Spec A$ is second countable iff $A$ is countable.
\item[(iii)]
$\Spec A$ is $T_0$, compact, and every closed irreducible set is the closure of a point.
\end{itemize}
\end{theorem}
\begin{proof}
The key point is that every filter $\ffrak$ is the intersection of all prime filters $\supseteq\ffrak$; this fact follows from a standard application of the Zorn Lemma. As a consequence, for every $D\subseteq A$, the intersection of all filters containing $D$ coincides with the intersection of all prime filters containing $D$.
This intersection, namely $\bigcap F_D$, is the smallest filter containing $D$, and we denote it by $\ffrak(D)$.
One verifies easily that $\ffrak(D)$ is the set of all $a\in A$ such that there exist $\vect a1r\in D$ satisfying
$a\ge a_1\tnorm\cdots\tnorm a_r$. Consider the mappings
\begin{align*}
A\supseteq D &\longmapsto F_D\in\text{sets closed in $\Spec A$} \\
\text{filters of $A$}\ni\textstyle{\bigcap} P
&\longleftarrow\!\mapstochar P\subseteq\Spec A
\end{align*}
They both reverse the $\subseteq$ relation. Their composition gives, on the left side, the mapping $D\mapsto \bigcap F_D=\ffrak(D)$ that associates to a set the filter it generates, and on the right side the topological closure mapping $P\mapsto F_{\bigcap P}$
(Proof: $\pfrak$ belongs to the topological closure of $P$ iff
$\forall a(\pfrak\in O_a\To O_a\cap P\not=\emptyset)$ iff
$\forall a(P\subseteq F_a\To a\in\pfrak)$ iff
$\pfrak\supseteq\bigcap P$).
As a consequence, they induce an antiisomorphism between the lattice of filters of $A$ and the lattice of closed sets of $\Spec A$; this proves~(i).

We leave~(ii) as an exercise, and prove~(iii). The $T_0$ property is clear, because the closure of $\pfrak$ is $F_{\bigcap\{\pfrak\}}=F_\pfrak=\{\qfrak:\qfrak\supseteq\pfrak\}$.
Compactness follows from~(ii) and the fact that $\Spec A=O_0$.
Let $F_\ffrak$ be a closed irreducible set, i.e., a closed set that cannot be expressed nontrivially as the union of two closed sets; we must show that $\ffrak$ is prime, i.e., that $a\lor b\in\ffrak$ implies ($a\in\ffrak$ or $b\in\ffrak$). Assume $a\lor b\in\ffrak$; then $F_{a\lor b}\supseteq F_\ffrak$. Since $F_{a\lor b}=F_a\cup F_b$, we have $F_\ffrak=(F_\ffrak\cap F_a)\cup(F_\ffrak\cap F_b)$, which must be a trivial decomposition. Hence either $F_a\supseteq F_\ffrak$ or
$F_b\supseteq F_\ffrak$, i.e., either $a\in\ffrak$ or $b\in\ffrak$.
\end{proof}

A \newword{spectral space} is a topological space in which the compact open sets form a basis closed under finite intersections, and such that the conditions in Theorem~\ref{ref2}(iii) hold. By~\cite{Hochster69}, these are exactly the prime ideal spaces of commutative rings with~$1$.

The spectral spaces of boolean algebras have a further property: the compact open sets are exactly the clopen sets (i.e., the sets which are both closed and open). Indeed, as we observed after Lemma~\ref{ref7}, if $A\in\Boole$ then every $\pfrak\in\Spec A$ is of the form $\pfrak=\varphi^{-1}[1]$ for some homomorphism $\varphi:A\to\two$. This immediately implies that $a\in\pfrak$ iff $\neg a\notin\pfrak$, i.e., $F_a=O_{\neg a}$. Therefore the $O_a$'s are clopen, and since every clopen is compact, there are no other clopens. 
Thus $A$ is in bijection with the clopen sets of $X=\Spec A$ via $a\mapsto F_a$, and since $F_{a\land b}=F_a\cap F_b$, $F_{\neg a}=X\setminus F_a$, $F_0=\emptyset$, and $F_1=X$, this bijection is a boolean algebra isomorphism.
We have thus proved the \newword{Stone Representation Theorem} \cite[\S18]{Halmos63}: every boolean algebra is isomorphic to the algebra of clopen subsets of its spectrum.

\begin{addendum}
Spectral\label{ref11} spaces can be functorially introduced in any congruence-modular equational class. In general, filters should be substituted by congruences (unless the class turns out to be
ideal-determined~\cite{GummUrsini84}), and one introduces the notion of prime congruence by using the commutator product; the construction carries on smoothly~\cite{agliano93}. The main trouble is with Theorem~\ref{ref2}(i): the open sets of $\Spec A$ will now be in 1-1 correspondence only with the radical congruences of $A$, i.e., those congruences $\Phi$ such that, for every congruence $\Psi$, if $\Phi$ contains the commutator product of $\Psi$ with itself, then $\Phi$ already contains $\Psi$ (the spectrum of $\Zbb$ as a commutative ring
is a typical example, the radical ideals being those generated by a squarefree integer). In our case there are no such problems, since equational classes are generated by truth-value algebras, in which a lattice structure is term-definable. All our classes are therefore congruence-distributive, and all congruences are radical.
\end{addendum}

\section{Proofs and dynamics}

We can now make precise the heuristic in the Introduction about the dynamics in Frege proof systems.

A \newword{substitution} on $FORM_\kappa$ is any mapping $\sigma:FORM_\kappa\to FORM_\kappa$ which distributes over the connectives (i.e., $\sigma(0)=0$, $\sigma(1)=1$, and $\sigma(r\circ s)=\sigma(r)\circ \sigma(s)$, for $\circ \in\{\tnorm,\to\}$). A substitution is therefore determined by an arbitrary assignment of formulas to propositional variables. Note that all variables must be substituted at the same time: e.g., if $r=x_1\to x_2$, $\sigma(x_1)=x_3\tnorm x_2$, and $\sigma(x_2)=x_1$, then $\sigma(r)=(x_3\tnorm x_2)\to x_1$. 

Given a formula $r$ and a set of formulas $\Theta$, a \newword{deduction} of $r$ from $\Theta$ is a finite sequence of formulas $\vect r1h$ such that $r_h=r$ and for every $1\le j\le h$ we have:
\begin{itemize}
\item[(a)] either $r_j\in\Theta$;
\item[(b)] or there exist $1\le k,m<j$ such that $r_m$ has the form $r_k\to r_j$ (we then say that $r_j$ follows from $r_k$ and $r_k\to r_j$ via Modus Ponens);
\item[(c)] or there exists $1\le k<j$ and a substitution $\sigma$ such that $r_j=\sigma(r_k)$.
\end{itemize}
An \newword{MP-deduction} is a deduction in which the substitution rule~(c) is never applied. For a fixed truth-value algebra $M$, it is often possible ---although sometimes difficult--- to find 
effectively a set of formulas $\Theta$ such that the formulas deducible from $\Theta$ are exactly the formulas which are true in $M$. If this happens, then we say that $\Theta$ provides an \newword{axiomatization} of $\operator{V}M$.

\begin{example}
Let $\Theta$ be the set of the following eight formulas:
\begin{gather*}
\bigl((x_0\to x_1)\tnorm(x_1\to x_2)\bigr)\to(x_0\to x_2); \\
(x_0\tnorm x_1)\to x_0; \quad 
0\to x_0; \\
(x_0\tnorm x_1)\to(x_1\tnorm x_0); \quad
(x_0\land x_1)\to(x_1\land x_0); \\
\bigl(x_0\to(x_1\to x_2)\bigr)\to\bigl((x_0\tnorm x_1)\to x_2\bigr);
\quad
\bigl((x_0\tnorm x_1)\to x_2\bigr)\to\bigl(x_0\to(x_1\to x_2)\bigr); \\
\bigl[\bigl((x_0\to x_1)\to x_2\bigr)\tnorm
\bigl((x_1\to x_0)\to x_2\bigr)\bigr]\to x_2.
\end{gather*}
All of them ---except perhaps the last one--- have rather transparent meanings: the first one expresses transitivity of implication, the third is \emph{ex falso quodlibet}, the fourth expresses commutativity of conjunction, and so on. We have~\cite{hajek98}:
\begin{itemize}
\item $\Theta\cup\{\neg\neg x_0\to x_0\}$ axiomatizes $\MV$;
\item $\Theta\cup\{\neg\neg x_0\to\bigl((x_1\tnorm x_0\to
x_2\tnorm x_0)\to (x_1\to x_2)\bigr), \neg(x_0\land\neg x_0)\}$ axiomatizes the equational class generated by
$M=\oi$ endowed with the product connectives in Example~\ref{ref3}(2);
\item $\Theta\cup\{x_0\to(x_0\tnorm x_0)\}$ axiomatizes the equational class generated by $M=\oi$ endowed with the G\"odel-Dummett connectives;
\item $\Theta\cup\{x_0\lor\neg x_0\}$ axiomatizes $\Boole$.
\end{itemize}
\end{example}

Fix a truth-value algebra $M$ and a cardinal $\kappa$. 
Let $\Theta_\kappa$ be the set of all formulas in $FORM_\kappa$ which are true in $M$. Hence $r(x_{i_1},\ldots,x_{i_n})\in\Theta_\kappa$ iff $r(a_{i_1},\ldots,a_{i_n})=1$ in $\free\kappa{\operator{V}M}$, where the $a_i$'s are the free generators. It is a customary abuse of notation to write $x_i$ for $a_i$, so the symbol $r$ may denote either an element of $FORM_\kappa$ or an element of $\free\kappa{\operator{V}M}$. This slight ambiguity is really sought for: it is exactly the ambiguity that results in working modulo $\Theta_\kappa$, or in identifying a formula with the function it induces on truth-values. Formally stated: $r,s\in FORM_\kappa$ are equal in $\free\kappa{\operator{V}M}$ iff $M\models r=s$ iff $r\to s,s\to r\in\Theta_\kappa$.

The set $\Theta_\omega$ is closed under Modus Ponens and substitution, so nothing new can be deduced from it. We let now $\Delta$ be any set of formulas, and raise two questions:
\begin{enumerate}
\item What can be deduced from $\Theta_\omega\cup\Delta$?
\item Which substitutions are needed for such a deduction?
\end{enumerate}

\begin{lemma}
If\label{ref8} $r$ is deducible from $\Theta_\omega\cup\Delta$, then there is a deduction involving only variables already appearing in $\{r\}\cup\Delta$.
\end{lemma}
\begin{proof}
By renaming variables we may assume that the variables appearing in $\{r\}\cup\Delta$ are exactly those variables with index $<\kappa$, for a certain $\kappa$.
Let a deduction of $r$ from $\Theta_\omega\cup\Delta$ be given. By a standard argument~\cite[p.~149]{church}, we may transform the given deduction into an MP-deduction $\vect r1h=r$ of $r$ from $\Theta_\omega\cup\{\sigma(t):\sigma\text{ is a substitution and }t\in\Delta\}$. Let $\tau$ be the substitution given by $\tau(x_j)=x_j$ if $j<\kappa$, and $\tau(x_j)=1$ otherwise. Then $\tau(r_1),\ldots,\tau(r_h)=r$ is an MP-deduction of $r$ from $\Theta_\kappa\cup\{\sigma(t):\sigma\text{ is a substitution, }t\in\Delta,\text{ and }\sigma(t)\in FORM_\kappa\}$, and hence a deduction of $r$ from $\Theta_\kappa\cup\Delta$ in which all formulas and all substitutions involve only variables with index $<\kappa$.
\end{proof}

Identifying $FORM_\kappa$ with $\free\kappa{\operator{V}M}$, a substitution is nothing more than an endomorphism of $\free \kappa{\operator{V}M}$ (i.e., a homomorphism from $\free \kappa{\operator{V}M}$ to itself). The freeness of the generators says that however we choose elements $\{b_i:i<\kappa\}$ in $\free \kappa{\operator{V}M}$ there is precisely one endomorphism that maps $x_i$ to $b_i$.

\begin{definition}
We denote by $\Sigma_\kappa$ the monoid of all endomorphisms of $\free\kappa{\operator{V}M}$, and by $\Xi_\kappa\subseteq\Sigma_\kappa$ the group of all automorphisms of $\free\kappa{\operator{V}M}$ (i.e., of the invertible elements of $\Sigma_\kappa$).
As explained before Theorem~\ref{ref2}, to every $\sigma\in\Sigma_\kappa$ 
there corresponds a continuous selfmapping $\sigma^*$ of $X_\kappa=\Spec\free\kappa{\operator{V}M}$, which we call the \newword{dual} of $\sigma$.
Let $\Pi\subseteq\Sigma_\kappa$, and let 
$O$ be an open subset of $X_\kappa$.
We define $(\Pi,O)$ to be the union of all backwards translates of $O$ under iteration of the substitutions in $\Pi$. Explicitly stated,
$$
(\Pi,O)=\bigcup\{(\sigma^*)^{-1}[O]:\sigma\text{ is in the submonoid of }\Sigma_\kappa\text{ generated by }\Pi\}.
$$
If $(\Pi,O)=X_\kappa$ for every $O\not=\emptyset$, then we say that $\Pi$ acts \newword{minimally} on $X_\kappa$: this is equivalent to saying that every point of $X_\kappa$ has a dense orbit under~$\Pi$.
\end{definition}

\begin{lemma}
Let\label{ref9} $\{r\}\cup\Delta\subseteq FORM_\kappa$. Then:
\begin{itemize}
\item[(i)] $r$ can be MP-deduced from $\Theta_\omega\cup\Delta$ iff $O_r\subseteq O_\Delta$ in $X_\kappa$;
\item[(ii)] if $\Delta'\subseteq FORM_\kappa$ is another set of formulas, then $O_\Delta=O_{\Delta'}$ iff $\Theta_\omega\cup\Delta$ and $\Theta_\omega\cup\Delta'$
MP-deduce the same formulas;
\item[(iii)] $r$ can be deduced from $\Theta_\omega\cup\Delta$ iff $O_r\subseteq(\Sigma_\kappa,O_\Delta)$.
\end{itemize}
\end{lemma}
\begin{proof}
(i) follows from Theorem~\ref{ref2}(i) and the proof of Lemma~\ref{ref8}: $r$ can be MP-deduced from 
$\Theta_\omega\cup\Delta$ iff $r$ can be MP-deduced from 
$\Theta_\kappa\cup\Delta$ iff $r$ belongs to the filter $\ffrak(\Delta)$ 
in $\free\kappa{\operator{V}M}$ iff $\ffrak(r)\subseteq\ffrak(\Delta)$ 
iff $O_r=O_{\ffrak(r)}\subseteq O_{\ffrak(\Delta)}=O_\Delta$. (ii) is 
clear: the sets $\Theta_\omega\cup\Delta$ and $\Theta_\omega\cup\Delta'$ MP-deduce 
the same formulas iff $\ffrak(\Delta)=\ffrak(\Delta')$ iff 
$O_\Delta=O_{\ffrak(\Delta)}=O_{\ffrak(\Delta')}=O_{\Delta'}$. We 
prove~(iii): assume that $r$ can be deduced from 
$\Theta_\omega\cup\Delta$. 
By Lemma~\ref{ref8}, there exists a deduction $\vect r1h=r$ 
such that 
every $r_j$ is in $\Theta_\kappa\cup\Delta$ and all substitutions 
applied are in $\Sigma_\kappa$. Working by induction on $h$ we assume 
that
$$
O_{r_1}\cup\cdots\cup O_{r_{h-1}}\subseteq (\Sigma_\kappa,O_\Delta).
$$
If $r_h\in\Delta$, then of course we are through (note that $s\in\Theta_\kappa$ is equivalent to $O_s=\emptyset$). For every $s,t\in\free\kappa{\operator{V}M}$ we have $O_t\subseteq O_s\cup O_{s\to t}$. Indeed, if $\pfrak\notin O_s$ and
$\pfrak\notin O_{s\to t}$, then $s,s\to t\in\pfrak$. Since filters are closed under MP, we have $t\in\pfrak$ and $\pfrak\notin O_t$. Therefore, if $r_h$ has been obtained via MP, then the induction hypothesis guarantees that $O_{r_h}\subseteq(\Sigma_\kappa,O_\Delta)$.
Finally, if $r_h=\sigma(r_j)$ for some $\sigma\in\Sigma_\kappa$ and $1\le j<h$, then $O_{r_h}=O_{\sigma(r_j)}=
(\sigma^*)^{-1}[O_{r_j}]\subseteq(\sigma^*)^{-1}[(\Sigma_\kappa,
O_\Delta)]\subseteq (\Sigma_\kappa,O_\Delta)$.
We leave the reverse implication as an exercise for the reader
(Hint: $O_r$ is compact).
\end{proof}

We say that $\Theta\subseteq FORM_\omega$ is \newword{equationally complete} if, however we choose $r,s\in FORM_\omega$ with $r\notin\Theta$, the formula $s$ is deducible from $\Theta\cup\{r\}$.

\begin{theorem}
The\label{ref10} set $\Theta_\omega$ is equationally complete iff $\Sigma_\omega$ acts minimally on $X_\omega$.
\end{theorem}
\begin{proof}
Let $O$ be a nonvoid open subset of $X_\omega$,
and assume that $\Theta_\omega$ is equationally complete: we want to show that $(\Sigma_\omega,O)=X_\omega$. Let $r\in\free\omega{\operator{V}M}$ be such that $\emptyset\not=O_r\subseteq O$; we then have $r\notin\Theta_\omega$. By assumption, the formula $0$ is deducible from $\Theta_\omega\cup\{r\}$, and hence we get $X_\omega=O_0\subseteq(\Sigma_\omega,O_r)\subseteq(\Sigma_\omega,O)$
from Lemma~\ref{ref9}(iii). The same argument yields the reverse implication.
\end{proof}

\begin{addendum}
The equational completeness of $\Theta_\omega$ amounts to the lack of nontrivial equational subclasses of $\operator{V}M$. Among classes generated by truth-value algebras, the only one fulfilling this property is $\Boole$ (easy proof, resting on the fact that $\two$ is a subalgebra of any $M$). Theorem~\ref{ref10} then implies that $\Sigma_\omega$ acts minimally on $X_\omega$ only in the case of boolean algebras.
There are other cases in which an algebra $M$ (not a truth-value algebra in our sense) generates a congruence-distributive equational class having no nontrivial subclasses. A particularly interesting case is when $M$ is the set of integers equipped with its natural structure $(\Zbb,+,-,0,\lor,\land)$ of lattice-ordered group~\cite{bkw}, \cite{andersonfei}. The resulting equational class is the class of all lattice-ordered groups, and the above properties are fulfilled. Theorem~\ref{ref10} says then that the endomorphisms of the free lattice-ordered groups act minimally on the relative spectra: see~\cite{pantiprime} for a description of such spaces.
\end{addendum}

For every
$p=(\ldots,p_i,\ldots)\in M^\kappa$, the evaluation mapping at $p$, given by 
$r(x_{i_1},\ldots,x_{i_n})\mapsto
r(p_{i_1},\ldots,p_{i_n})$, is a homomorphisms
$\varphi:\free\kappa{\operator{V}M}\to M$. Since $M$ is totally-ordered, the kernel $\varphi^{-1}[1]$ is a prime filter $\pfrak$, hence an element of $X_\kappa$. We thus get a mapping $p\mapsto\pfrak$ from $M^\kappa$ to $X_\kappa$, which we denote by $\pi$.
The map $\pi$
has dense range ($\emptyset\not= O_r$ $\To$ $r\notin\Theta_\kappa$ $\To$ $\exists p\in
M^\kappa\;r(p)\not=1$ $\To$ $\exists p\; r\notin\pi(p)$
$\To$ $\exists p\;\pi(p)\in O_r$), but it
is not necessarily continuous; it is continuous in the two cases that most concern us, namely in classical logic (see Lemma~\ref{ref25}) and in \Luk\ logic (see the next section).

Let $\sigma:\free\kappa{\operator{V}M}\to\free\kappa{\operator{V}M}$ be a substitution, $s_i=\sigma(x_i)$. Then the $\kappa$-tuple $(\ldots,s_i,\ldots)$ determines a function $S:M^\kappa\to M^\kappa$ via $p\mapsto(\ldots,s_i(p),\ldots)$. Since $\pi(S(p))=\{r\in\free\kappa{\operator{V}M}:
r(\ldots,s_i(p),\ldots)=1\}=\{r:[\sigma(r)](p)=1\}=
\sigma^{-1}(\{t:t(p)=1\})=\sigma^*(\pi(p))$, the diagram
\begin{equation}\tag{$*$}
\begin{split}
\begin{xy}
\xymatrix{
{M^\kappa} \ar[r]^S \ar[d]_{\pi} & {M^\kappa} \ar[d]^{\pi} \\
{X_\kappa} \ar[r]_{\sigma^*} & {X_\kappa}
}
\end{xy}
\end{split}
\end{equation}
commutes. As $\pi$ has dense range, $S$ determines the continuous function~$\sigma^*$.
We call $\sigma^*$ the \newword{dual} of $\sigma$, and $S$ the \newword{mapping on truth-values} induced by $\sigma$.

\begin{lemma}
In\label{ref25} classical logic, the map $\pi:\two^\kappa
\to X_\kappa=\Spec\free\kappa\Boole$ is a homeomorphism.
\end{lemma}
\begin{proof}
Recall from the end of Section~\ref{ref24} that $\free\kappa\Boole$ is the boolean algebra of all clopen subsets of $\two^\kappa$, where
the latter space is given the product topology. If $p$ and $q$ are distinct points of $\two^\kappa$, then there exists a clopen set $C$ containing $p$ and not containing $q$. Hence $C\in\pi(p)\setminus\pi(q)$ and $\pi$ is injective (as usual, we are identifying clopen subsets with their characteristic functions). Let $\pfrak\in X_\kappa$. Since $\pfrak$ is a proper filter, $\emptyset\notin\pfrak$ and in particular the intersection
of any finite family of elements of $\pfrak$ is nonempty. By compactness, $\bigcap\pfrak$ contains a point $p\in\two^\kappa$.
Since $C\in\pfrak$ implies $p\in C$, we have $\pfrak\subseteq\pi(p)$.
But in a boolean algebra every prime filter is maximal, and hence 
$\pfrak=\pi(p)$. So $\pi$ is a bijection. Both in $\two^\kappa$
and in $X_\kappa$ the clopen sets generate the topology; moreover, as shown before Addendum~\ref{ref11}, the mapping $C\mapsto F_C$ is a bijection between the two families of clopen sets. Since $p\in C$ iff $C\in\pi(p)$ iff $\pi(p)\in F_C$, the map $\pi$ is a homeomorphism.
\end{proof}

Given a point $p$ and a nonvoid open subset $O$ of the Cantor space
$\two^\omega$, one easily constructs a homeomorphism $S:\two^\omega\to\two^\omega$ such that $S(p)\in O$ (if $p=(p_0,p_1,\ldots)$ and $[\vect a0n]=\{q\in\two^\omega:q_i=a_i\text{ for }i=0,\ldots,n\}$ is a block contained in $O$, then the mapping that exchanges $0$ with $1$ in those indices $i$ for which $p_i\not=a_i$ is such a homeomorphism).
Hence not only $\Sigma_\omega$, but even $\Xi_\omega$ acts minimally on $\Spec\free\omega\Boole$.

Of course we can do better than that, because there exist many minimal homeomorphisms of the Cantor space. The simplest example is obtained by identifying $\two^\omega$ with the topological group of $2$-adic integers $\Zbb_2$, and letting $S$ be the translation by~$1$: $S(p)=p+1$. Let us compute the substitution $\sigma$ on $\free\omega\Boole$ for which $S=\sigma^*$. If $x_i$ is the $i$-th free generator, then $F_{x_i}=\{p\in\two^\omega:p_i=1\}$, and $F_{\sigma(x_i)}=(\sigma^*)^{-1}[F_{x_i}]=\{p:S(p)\in F_{x_i}\}=\{p:(p+1)_i=1\}$. Since addition in $\Zbb_2$ is just addition in base $2$ with carry, we have that $p\in F_{\sigma(x_i)}$ iff
\begin{itemize}
\item either $p_i=1$ and $p_j=0$ for some $j<i$;
\item or $p_i=0$ and $p_j=1$ for every $j<i$.
\end{itemize}
Therefore
\begin{align*}
F_{\sigma(x_i)} &=
\bigl[F_{x_i}\cap(O_{x_0}\cup\cdots\cup O_{x_{i-1}})\bigr] \\
&\quad \cup\bigl[O_{x_i}\cap F_{x_0}\cap\cdots\cap F_{x_{i-1}}\bigr] \\
&= \bigl[F_{x_i}\cap(F_{\neg x_0}\cup\cdots\cup
F_{\neg x_{i-1}})\bigr] \\
&\quad \cup\bigl[F_{\neg x_i}\cap F_{x_0}\cap\cdots\cap
F_{x_{i-1}}\bigr].
\end{align*}
Consider the following formulas:
\begin{align*}
s_0 &= \neg x_0 \\
s_i &= \bigl[x_i\land(\neg x_0\lor\cdots\lor\neg x_{i-1})\bigr]
\lor\bigl[\neg x_i\land x_0\land\cdots\land x_{i-1}\bigr] \\
&= x_i\triangle(x_0\land\cdots\land x_{i-1}) 
\quad\text{(for $i>0$)}
\end{align*}
($\triangle$ is the boolean symmetric difference: $a\triangle b=(a\land\neg b)\lor(\neg a\land b)$). Then $F_{\sigma(x_i)}=F_{s_i}$ by the isomorphism cited before Addendum~\ref{ref11}. The required substitution is therefore the one defined by $\sigma(x_i)=s_i$. From the point of view of proof systems, we have thus obtained the following result.

\begin{theorem}
From\label{ref14} the set of boolean tautologies plus any given non-tautology we can derive every formula using only Modus Ponens and the substitution $\sigma$ given above.
\end{theorem}

\section{\Luk\ logic}

In the rest of this paper we will concentrate on \Luk\ logic; we therefore fix $M=\oi$ endowed with the connectives in Example~\ref{ref3}(3).
A key distinguishing feature of the \Luk\ connectives is their continuity with respect to the standard topology of $\oi$. As a matter of fact \Luk\ logic is the only t-norm based logic in which all connectives are continuous~\cite{menupav}.

A \newword{[rational] cellular complex over $M^n=\oi^n$} is a finite set $W$ of \newword{cells} (i.e., compact convex polyhedrons), whose union is $\oi^n$, and such that:
\begin{enumerate}
\item every vertex of every cell of $W$ has rational coordinates;
\item if $C\in W$ and $D$ is a face of $C$, then $D\in W$;
\item every two cells intersect in a common face.
\end{enumerate}
A \newword{McNaughton function} is a continuous function $f:\oi^n\to\oi$ for which there exists a complex 
as above and affine linear functions with integer coefficients $F_j(\bar x)=a_j^1x_1+\cdots+a_j^nx_n+a_j^{n+1}$, in 1-1 correspondence with the $n$-dimensional cells $C_j$ of the complex, such that $f\restriction C_j=F_j$ for each $j$.

\begin{theorem}
\textrm{\cite{mcnaughton}, \cite{mundicimn}}
The\label{ref15} elements of $\free n\MV$ (i.e., the functions from $\oi^n$ to $\oi$ induced by a formula of \Luk\ logic) are exactly the McNaughton functions.
\end{theorem}

Here are typical McNaughton functions, for $n=1$ and $n=2$: they are induced by the formulas $\neg x_0\lor\bigl((x_0\land\neg x_0)\oplus
(x_0\land\neg x_0)\bigr)$ and $(x_0\to x_1)\land(x_0\oplus x_0\oplus x_1\oplus x_1)$, respectively.
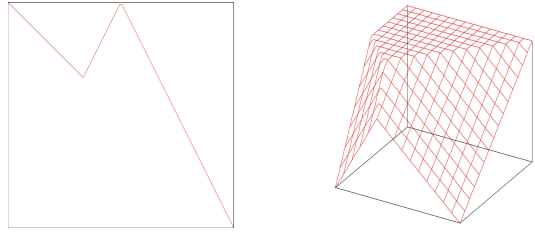
\begin{figure}[H]
\begin{center}
\includegraphics[width=3cm,height=3cm]{figura1.epsf}
\quad\quad\quad
\includegraphics[width=3cm,height=3cm]{figura2clip.epsf}
\end{center}
\end{figure}

Given a substitution $\sigma:\free n\MV\to\free n\MV$, to each function $s_i=\sigma(x_i)$ there corresponds a cellular complex $W_i$ such that $s_i$ is affine linear on each cell of~$W_i$. Let $W$ be a complex that is a common refinement of $\vect W0{{n-1}}$. Then on each cell $C_j$ of $W$ the function $S:\oi^n\to\oi^n$ defined before Lemma~\ref{ref25} is given by
\begin{equation}\tag{$**$}
\begin{pmatrix}
p_0 \\
\vdots \\
p_{n-1}
\end{pmatrix}
\mapsto
A_j\begin{pmatrix}
p_0 \\
\vdots \\
p_{n-1}
\end{pmatrix}
+B_j,
\end{equation}
where $A_j$ is an $n\times n$ matrix and $B_j$ a column vector, both having integer coefficients.
Conversely, every continuous selfmapping $S$ of $\oi^n$ which is piecewise affine linear with integer coefficients (i.e., is locally expressible in the form $(**)$, using finitely many $A_j$'s and $B_j$'s) is induced by some endomorphism $\sigma$ of $\free n\MV$.
We call such an $S$ a \newword{McNaughton mapping}; if moreover $S$ is invertible we call it a \newword{McNaughton homeomorphism}; McNaughton homeomorphisms on $\oi^n$ are exactly the mappings on truth-values induced by the automorphisms of
$\free n\MV$.
Apart from its relevance in \Luk\ logic, the class of McNaughton mappings is quite interesting \emph{per se}.

\begin{example}
Consider the following complexes over $\oi^2$; they are both
symmetric under a $\pi$ rotation about the centre of the square.
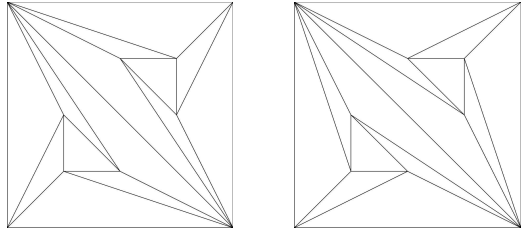
\begin{figure}[H]
\begin{center}
\includegraphics[height=3cm,width=3cm]{homeom-1.epsf}
\qquad
\includegraphics[height=3cm,width=3cm]{homeom-2.epsf}
\end{center}
\end{figure}
The vertices of the lower inner triangle are
$p_0=(1/4,1/4)$, $p_1=(1/2,1/4)$, $p_2=(1/4,1/2)$; for $0\le i\le 2$, let $p_i'$ be the vertex symmetric to $p_i$.
Then there exists a unique homeomorphism $S$ such that:
\begin{enumerate}
\item $S(p_i)=p_{i+1\pmod{3}}$, and $S(p'_i)=p'_{i+1\pmod{3}}$;
\item every other vertex is fixed;
\item $S$ is affine linear on each cell.
\end{enumerate}
In short, the first complex is mapped onto the second by ``rotating counterclockwise'' the two inner triangles, and distorting accordingly the border triangles. As a matter of fact, $S$ is topologically conjugate to the union of two  twists~\cite[\S5]{pantibernoulli}.
The data above determine the matrix $A_j$ and the column vector $B_j$ on each triangle $C_j$. One checks directly that all these matrices and vectors have integer entries; hence $S$ is a McNaughton homeomorphism.
In doing computations, it is expedient to write
$p=(\vect p0{{n-1}})\in\oi^n$ using projective coordinates
$(\projvect a0n)\sim(\projvect p0{{n-1}}:1)$. For example, if $C_1$ is the triangle $\angles{p_0,(1,0),p_1}$, which is mapped to
$\angles{p_1,(1,0),p_2}$, then $A_1$ and $B_1$ are the upper left $2\times 2$ matrix and upper right $2\times 1$ column vector in the matrix
$$
\begin{pmatrix}
-1 & -5 & 2 \\
1 & 4 & -1 \\
0 & 0 & 1
\end{pmatrix}
=
\begin{pmatrix}
2 & 1 & 1 \\
1 & 0 & 2 \\
4 & 1 & 4
\end{pmatrix}
\begin{pmatrix}
1 & 1 & 2 \\
1 & 0 & 1 \\
4 & 1 & 4
\end{pmatrix}^{-1}.
$$
\end{example}

If $S$ is induced by $\sigma$, then the action of $S$ on $\oi^n$ is just the surface of the action of the full dual $\sigma^*$ on $X_n$. Indeed, as proved in~\cite[Proposition~8.1]{mundicijfa},
in the case of \Luk\ logic
the map $\pi$ in the diagram $(*)$ is a homeomorphic embedding of $\oi^n$ onto the subspace of maximal filters.

By Lemma~\ref{ref7}, the points of $X_n$ (indeed, of any spectrum) form a forest under the \newword{specialization order}: $\pfrak\le\qfrak$ iff $\qfrak$ is in the closure of $\{\pfrak\}$ iff
$\pfrak\subseteq\qfrak$ (a \newword{forest}, sometimes called a \newword{root system}, is a poset in which the elements greater than any given element form a chain).

\begin{example}
A full description of $X_n$ is given
in~\cite{pantiprime}. Since it is rather involved, here we limit ourselves to the cases $n=1$ and $n=2$.

Assume $n=1$. If $p\in(0,1)$ is rational, then there are two incomparable prime filters $\pi(p)^+$ and $\pi(p)^-$ properly contained in the maximal $\pi(p)$. Namely, $\pi(p)^-$ is the filter of all McNaughton functions $:\oi\to\oi$ that are $1$ in a left neighborhood of $p$, and analogously for $\pi(p)^+$ w.r.t.~right neighborhoods. The only prime filter contained in $\pi(0)$ (respectively, in $\pi(1)$) is $\pi(0)^+$ (respectively, $\pi(1)^-$). If $p$ is irrational, then $\pi(p)$ is minimal in the specialization order.

Now assume $n=2$, $p=(p_0,p_1)\in(0,1)^2$. If $p_0,p_1,1$ are linearly independent over $\Qbb$, then $\pi(p)$ is a minimal ---as well as maximal--- prime filter. If $p_0,p_1,1$ satisfy exactly one (up to scalar multiples) nontrivial linear dependence over $\Qbb$, then there are two incomparable prime filters below $\pi(p)$, and both of them are minimal. Otherwise, consider the unit circle $S^1$ in the tangent space to $p$. For every $u\in S^1$ there is a prime filter $\pfrak_u$ contained in $\pi(p)$. If the line in the tangent space connecting the origin with $u$ does not hit any point with rational coordinates, then $\pfrak_u$ is minimal. Otherwise, $\pfrak_u$ contains two minimal prime filters $\pfrak_u^+$ and $\pfrak_u^-$.
For $p$ along the border of the unit square, this description gets modified in the obvious way. The following picture may clarify the situation:
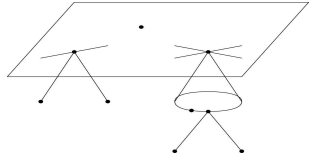
\begin{figure}[H]
\begin{center}
\includegraphics[width=4cm,height=2cm]{figura-1.epsf}
\end{center}
\end{figure}
\end{example}

\begin{addendum}
Let $n$ (not necessarily $n=1,2$) be given, $\pfrak\in X_n$. Let $\pfrak=\pfrak_0\subset \pfrak_1\subset\cdots
\subset\pfrak_t$ be the chain, of length $t$, of elements above $\pfrak$ in the specialization order.
Given an endomorphism $\sigma$ of $\free n\MV$, we have
$$
\frac{\free n\MV}{\sigma^{-1}[\pfrak]} \simeq
\frac{\sigma[\free n\MV]}{\pfrak\cap\sigma[\free n\MV]} \subseteq
\frac{\free n\MV}{\pfrak},
$$
and hence the MV-algebra $\free n\MV/\sigma^*(\pfrak)$ is a subalgebra of $\free n\MV/\pfrak$. By~\cite[Theorem~4.7(i) and Corollary~4.9]{pantiprime}, this implies that the length of the chain above $\sigma^*(\pfrak)$ is less than or equal to $t$.
\end{addendum}

Since $\oi^n$ is dense in $X_n$, in principle we might reduce the study of $\sigma^*$ to the study of $S$. However, taking into consideration the action of $\sigma^*$ on the full spectrum gives us deeper insight. For example, it is possible to provide an intrinsic (i.e., coordinate-free) characterization of the differentials $T_pS$ of a McNaughton mapping $S$. Differentials of piecewise-linear maps have been constructed by Tsujii in~\cite{Tsujii01}; we show here how Tsujii's construction can be intrinsically described in purely algebraic terms.

It may be helpful for the reader to recall the coordinate-free description of the differentials of a morphism $S:X\to Y$ of differentiable varieties. Let $p$ be a point of $X$, $q=S(p)$, $\Ocal_p$ and $\Ocal_q$ the rings of germs of differentiable functions at $p$ and $q$, respectively. $\Ocal_p$ and $\Ocal_q$ are local rings: let $\mfrak$ and $\nfrak$ be the respective maximal ideals (i.e., $\mfrak=\{f\in\Ocal_p:f(p)=0\}$, and analogously for $\nfrak$).
The mapping $\sigma:\Ocal_q\to\Ocal_p$ defined by $\sigma(g)=g\circ S$ is a well-defined ring homomorphism, and $\sigma[\nfrak]\subseteq\mfrak$. Therefore, $\sigma$ induces a vector space homomorphism from $\nfrak/\nfrak^2$ to $\mfrak/\mfrak^2$, which we denote by $\bar{\sigma}$. The tangent spaces $T_pX$ and $T_qY$ are canonically isomorphic to the dual vector spaces $(\mfrak/\mfrak^2)'$
and $(\nfrak/\nfrak^2)'$, respectively, and under these isomorphisms the differential $T_pS$ corresponds to the dual mapping 
$\bar{\sigma}':(\mfrak/\mfrak^2)'\to(\nfrak/\nfrak^2)'$. 
Explicitly, if
$T_pX\ni v:\mfrak/\mfrak^2\to\Rbb$ is a tangent vector at $p$, then
$(T_pS)(v)$ is the tangent vector at $q$ defined by
$[(T_pS)(v)](g/\nfrak^2)=(v\circ\bar{\sigma})(g/\nfrak^2)=
v((g\circ S)/\mfrak^2)$.

We will develop an analogous description for piecewise-linear maps.
Before doing so, we need a few more preliminaries; see~\cite{bkw}, \cite{andersonfei} for more details and unproved claims. A \newword{lattice-ordered abelian group} (\newword{\llgroup} for short) is a structure
$(G,+,-,0,\land,\lor)$ such that $(G,+,-,0)$ is an abelian group,
$(G,\land,\lor)$ is a lattice, and $+$ distributes over the lattice operations. \llhomomorphisms\ of \llgroups\ are groups homomorphisms that are also lattice homomorphisms. The class of all \llgroups\ is equational and, as such, contains free objects. The free \llgroup\ over $n$ generators, $\Fl(n)$, is the \llgroup\ (under pointwise operations) of all functions $g:\Rbb^n\to\Rbb$ that are continuous and piecewise-linear with integer coefficients (i.e., there exist
finitely many homogeneous linear polynomials
$\vect g1m\in\Zbb[\vect x1n]$ such that, for every $w\in\Rbb^n$,
$g(w)=g_j(w)$ for some $1\le j\le m$). The set $\lHom(\Fl(n),\Rbb)$
of all \llhomomorphisms\ from $\Fl(n)$ to $\Rbb$ is in 1-1
correspondence with $\Rbb^n$, via the map that associates
to $w\in\Rbb^n$ the evaluation mapping $\varphi_w:g\mapsto g(w)$. A \newword{strong unit} of the \llgroup\ $G$ is an element $0\le u\in G$
such that, for every $g\in G$, $g\le nu$ for some positive integer $n$.
If $u$ is a strong unit of $G$, then the interval $[0,u]=
\{g\in G:0\le g\le u\}$ can be given the structure of an MV-algebra
$\Gamma(G,u)=([0,u],\oplus,\neg,0,1)$
by setting $g\oplus h=(g+h)\land u$, $\neg g=u-g$, $0=0_G$, 
$1=u$. The mapping $(G,u)\mapsto\Gamma(G,u)$ is functorial, and determines a categorical equivalence between the category of \llgroups\ with strong unit and the category of MV-algebras~\cite{mundicijfa}. In particular, the filters of $\Gamma(G,u)$ are in natural 1-1 correspondence with the kernels of
\llhomomorphisms\ of domain $G$.

The preliminaries being over, let $S:\oi^n\to\oi^n$ be a McNaughton mapping, $p\in\oi^n$, $q=S(p)$. Then $\mfrak=\pi(p)$ and $\nfrak=\pi(q)$ are maximal filters of $\free n\MV$. Given $\pfrak\in X_n$, the \newword{germinal filter} corresponding to $\pfrak$ is the filter
$\gfrak_\pfrak=\bigcap\{\qfrak\in X_n:\qfrak\subseteq\pfrak\}$.
By~\cite[Proposition~10.5.3 and Definition~10.5.6]{bkw}, and using the properties of the $\Gamma$ functor, the quotient $A_p=\free n\MV/\gfrak_\mfrak$ is the MV-algebra of germs at $p$ of McNaughton functions; analogously for
$A_q=\free n\MV/\gfrak_\nfrak$. The MV-algebra $A_p$ is \newword{local}, i.e., has a unique maximal filter $\mfrak/\gfrak_\mfrak$.
Let $\sigma$ be the endomorphism of $\free n\MV$ that induces $S$.

\begin{lemma}
Notation being as above, $\sigma[\nfrak]\subseteq\mfrak$ and $\sigma[\gfrak_\nfrak]\subseteq\gfrak_\mfrak$.
\end{lemma}
\begin{proof}
By the commutativity of the diagram $(*)$, $\sigma^{-1}[\mfrak]=
\sigma^*(\mfrak)=\nfrak$, so the first statement is immediate.
Let $\pfrak$ be a prime filter below $\mfrak$ in the specialization order. Since $\sigma^*$ is continuous and $\mfrak$ is in the closure of $\{\pfrak\}$, the maximal filter $\nfrak=\sigma^*(\mfrak)$ must be in the closure of $\{\sigma^*(\pfrak)\}$. Therefore $\sigma^*(\pfrak)$ is below $\nfrak$, hence $\gfrak_\nfrak\subseteq\sigma^*(\pfrak)=
\sigma^{-1}[\pfrak]$ and $\sigma[\gfrak_\nfrak]\subseteq\pfrak$.
\end{proof}

As a consequence, $\sigma$ determines a homomorphism of \llgroups
$$
\bar{\sigma}:\nfrak/\gfrak_\nfrak
\to \mfrak/\gfrak_\mfrak,
$$
which plays the r\^ole of the codifferential on cotangent spaces.
Since the composition of \llhomomorphisms\ is an \llhomomorphism, $\bar{\sigma}$ induces a dual mapping
$\bar{\sigma}':\lHom(\mfrak/\gfrak_\mfrak,\Rbb)\to
\lHom(\nfrak/\gfrak_\nfrak,\Rbb)$ by
$\bar{\sigma}'(\varphi)=\varphi\circ\bar{\sigma}$.

\begin{theorem}
Under\label{ref22} the identification $w\mapsto\varphi_w$ of
$\Rbb^n$ with $\lHom(\Fl(n),\Rbb)$
described above, the map $\bar{\sigma}'$ corresponds to Tsujii's
differential.
\end{theorem}
\begin{proof}
For simplicity's sake, we assume that $p$ and $q=S(p)$ have rational coordinates and are in the topological interior of the $n$-cube: we will discuss in Addendum~\ref{ref21} how these assumptions can be discarded.
Write $q$ in projective coordinates $(\projvect a0n)$, with
$a_n>0$, and let $Q=(\vect a0n)\in\Rbb^{n+1}$.
Let $\Nfrak$ be the kernel of $\varphi_Q$,
and let $\Gfrak_\Nfrak$ be the germinal
kernel associated to $\Nfrak$~\cite[Proposition~10.5.3]{bkw}.
Since $q$ has rational coordinates, $Q$ has rank $1$ according to~\cite[p.~188]{pantiprime}. By~\cite[Theorem~4.8]{pantiprime},
the quotient $\Nfrak/\Gfrak_\Nfrak$ is an \llgroup, which is \llisomorphic\ to $\Fl(n)$ under the map $D_Q$ defined in~\cite[Definition~2.2]{pantiprime}. By the
properties of the $\Gamma$ functor, the \llgroups\ $\nfrak/\gfrak_\nfrak$
and $\Nfrak/\Gfrak_\Nfrak$ are \llisomorphic\ as well. 
We therefore obtain an \llisomorphism\ $D_q:\nfrak/\gfrak_\nfrak\to
\Fl(n)$ which, by explicit computation, has the form
$$
\bigl[D_q(r/\gfrak_\nfrak)\bigr](w)=
\lim_{h\to0^+}\frac{r(q+hw)-r(q)}{h}.
$$
We have of course an analogous \llisomorphism\ 
$D_p:\mfrak/\gfrak_\mfrak\to\Fl(n)$.
Let $\mathcal{D}_pS$ denote the Tsujii differential of $S$ at $p$,
and let $v\in\Rbb^n$. Since $S$ is continuous and defined everywhere
on $\oi^n$, the definition in~\cite[Eq.~(13)]{Tsujii01} simplifies to
$$
(\mathcal{D}_pS)(v)=
\lim_{h\to0^+}\frac{S(p+hv)-S(p)}{h}.
$$
We want to show that $\bar{\sigma}'(\varphi_v\circ D_p)=
\varphi_{(\mathcal{D}_pS)(v)}\circ D_q$. Setting
$(\mathcal{D}_pS)(v)=w\in\Rbb^n$, this amounts to the commutativity
of the diagram
$$
\begin{xy}
\xymatrix{
{\nfrak/\gfrak_\nfrak} \ar[0,2]^{\bar{\sigma}} \ar[d]_{D_q} & & 
{\mfrak/\gfrak_\mfrak} \ar[d]^{D_p} \\
{\Fl(n)} \ar[r]_{\varphi_w} & {\Rbb} & 
{\Fl(n)} \ar[l]^{\varphi_v}
}
\end{xy}
$$
Choose $r/\gfrak_\nfrak\in\nfrak/\gfrak_\nfrak$. As remarked
in~\cite[p.~358]{Tsujii01}, the terms $h^{-1}\bigl(r(q+hw)-r(q)\bigr)$, $h^{-1}\bigl((r\circ S)(p+hv)-(r\circ S)(p)\bigr)$,
and $h^{-1}\bigl(S(p+hv)-S(p)\bigr)$ take constant values for
sufficiently small $h>0$. If $h_0$ is such an $h$, we get
\begin{align*}
\varphi_w\bigl(D_q(r/\gfrak_\nfrak)\bigr) &=
\bigl(D_q(r/\gfrak_\nfrak)\bigr)(w) \\
&= h_0^{-1}\bigl(r(q+h_0w)-r(q)\bigr) \\
&= h_0^{-1}\bigl(r(q+S(p+h_0v)-S(p))-r(q)\bigr) \\
&= h_0^{-1}\bigl((r\circ S)(p+h_0v)-(r\circ S)(p)\bigr) \\
&= \bigl(D_p((r\circ S)/\gfrak_\mfrak)\bigr)(v) \\
&= \varphi_v\bigl(D_p(\bar{\sigma}(r/\gfrak_\nfrak))\bigr),
\end{align*}
as required.
\end{proof}

\begin{addendum}
In\label{ref21} the proof of Theorem~\ref{ref22} we assumed that $p$ and $q$ have rational coordinates and are in the topological interior of $\oi^n$. The first assumption is motivated by the fact that, by definition, McNaughton mappings have integer coefficients.
This implies that the only tangent vectors, say at $q$, that can be algebraically recognized are those in the $\Rbb$-span $V$ of the set of
all $v\in\Rbb^n$ such that the affine line through $q$ and $q+v$
is the intersection of affine hyperspaces having integer (equivalently, rational) coefficients. The \llgroup\ $\nfrak/\gfrak_\nfrak$ is then \llisomorphic\ to $\Fl(k)$, for $k=\dim(V)\le n$, and the algebraic tangent space $(\nfrak/\gfrak_\nfrak)'$ is isomorphic to $\Rbb^k$. According to the personal interests, one may either accept these underdimensional tangent spaces, or treat them drastically, by tensoring everything
with $\Rbb$. This means considering piecewise-linear functions with arbitrary real coefficients, so dropping the assumption that the cellular complexes involved have rational vertices, and
passing from \llgroups\ and MV-algebras to real vector lattices~\cite{baker68} and their $\Gamma$ images. All quotients
$\nfrak/\gfrak_\nfrak$ are then isomorphic to the free vector lattice
over $n$ generators $\FVL(n)$~\cite[Theorem~3.8]{pantiprime}.
The dual of $\FVL(n)$, i.e., the set of all $\Rbb$-linear
\llhomomorphisms\ from $\FVL(n)$ to $\Rbb$ is still in bijection with $\Rbb^n$ via the evaluation mapping, and all dimensionality problems disappear.

About the other assumption: if $p$ or $q$ (say $p$) is on the boundary of $\oi^n$, then the quotient $\mfrak/\gfrak_\mfrak$ is \llisomorphic\
not to $\Fl(n)$, but to a quotient of $\Fl(n)$ by a principal
kernel or, equivalently, to the \llgroup\ of restrictions of the
elements of $\Fl(n)$ to a polyhedral cone $W$. The dual $(\mfrak/\gfrak_\mfrak)'$ is then in bijection with $W$, in agreement with~\cite[Eq.~(14)]{Tsujii01}, and the proof of Theorem~\ref{ref22} carries on.
\end{addendum}

\section{Chaotic actions}

Let $p=(\ldots,p_i,\ldots)\in\oi^\kappa$. We say that $p$ has \newword{finite denominator} if all $p_i$'s are rational numbers, and there exists $0<d\in\Zbb$ such that $dp\in\Zbb^\kappa$. The least such $d$ is the \newword{denominator} of $p$, written $\den(p)$.

\begin{lemma}
Let\label{ref16} $p,q\in\oi^\kappa$.
\begin{itemize}
\item[(i)] If $p$ has finite denominator and $\sigma\in\Sigma_\kappa$, then $S(p)$ has finite denominator and $\den(S(p))\mid\den(p)$.
\item[(ii)] If $\den(q)\mid\den(p)$, then $S(p)=q$ for some $\sigma\in\Sigma_\kappa$.
\item[(iii)] If $p$ does not have finite denominator, then the $\Sigma_\kappa$-orbit of $p$ is dense.
\end{itemize}
\end{lemma}
\begin{proof}
(i) Let $d=\den(p)$, $S(p)=(\ldots,q_i,\ldots)$, $\sigma(x_i)=
s_i(\vect x{{j_1}}{{j_n}})$. Then $q_i=s_i(\vect p{{j_1}}{{j_n}})=
a^1p_{j_1}+\cdots+a^np_{j_n}+a^{n+1}$, for some 
$a^1,\ldots,a^{n+1}\in\Zbb$.
Therefore $dq_i\in\Zbb$ and $dS(p)\in\Zbb^\kappa$, and our claim easily follows.\newline
(ii) Let $d=\den(p)$. The integers $\{dp_i:i<\kappa\}\cup\{d\}$ must be relatively prime (otherwise $\den(p)$ would be smaller than $d$).
Therefore there exist indices $\vect j1m$ and integer numbers
$a^1,\ldots,a^m,a^{m+1}$ such that $d(a^1p_{j_1}+\cdots+
a^mp_{j_m}+a^{m+1})=1$. Let $n\le\kappa$ be greater than $\vect j1m$. Then the affine linear polynomial
$f_i=q_id(a^1x_{j_1}+\cdots+a^mx_{j_m}+a^{m+1})$ has integer coefficients, since $\den(q)\mid d$. The function $(f_i\lor 0)\land 1:
\oi^n\to\oi$ is a McNaughton function, hence by Theorem~\ref{ref15}
it is expressible via a formula $s_i\in\free n\MV$. Since
$s_i(\vect p{{j_1}}{{j_n}})=q_i$, the substitution $\sigma\in\Sigma_\kappa$ defined by $\sigma(x_i)=s_i$ satisfies our requirements.\newline
(iii) It suffices to show that, for every $0\le a<b\le 1$, there exists an element $s\in\free\kappa\MV$ such that $a<s(p)<b$.
Let $G$ be the additive subgroup of $\Rbb$ generated by $\{p_i:i<\kappa\}\cup\{1\}$. Since $p$ does not have finite denominator, $G$ is dense in $\Rbb$, and therefore there exist indices $\vect j1m$
and integer numbers
$a^1,\ldots,a^m,a^{m+1}$ such that $a<a^1p_{j_1}+\cdots+
a^mp_{j_m}+a^{m+1}<b$. One then argues as in~(ii) above.
\end{proof}

In the following we will tacitly identify via $\pi$ the $\kappa$-cube $\oi^\kappa$ with the subspace of $X_\kappa$ whose elements are the maximal filters.

\begin{corollary}
Let\label{ref17} $\Rat_n$ be the set of rational points in $\oi^n$. Then $\Rat_n$ is a dense subset both of $\oi^n$ and of $X_n$. All the 
elements of $\Rat_n$ have a finite $\Sigma_n$-orbit. No point of $X_\omega$ has a finite $\Sigma_\omega$-orbit.
\end{corollary}
\begin{proof}
$\Rat_n$ coincides with the set of points in $\oi^n$ having finite denominator, and is dense in $\oi^n$. As the $n$-cube is dense in $X_n$, $\Rat_n$ is dense in $X_n$ as well. If $p\in\Rat_n$, then by Lemma~\ref{ref16} the $\Sigma_n$-orbit of $p$ is the set of points whose denominator divides $\den(p)$, and this set is finite. Since the submonoid of $\Sigma_\omega$ whose elements are all the substitutions $\sigma$ such that $\sigma(x_i)\in\ooii$ has the cardinality of the continuum, our last claim is immediate.
\end{proof}

By Lemma~\ref{ref16}(i) $\Sigma_\kappa$ does not act minimally on $X_\kappa$. We are therefore lead to weaken the requirement of minimality to that of topological transitivity: we say that $\Pi\subseteq\Sigma_\kappa$ is \newword{topologically transitive} on $X_\kappa$ if $(\Pi,O)$ is dense in $X_\kappa$, for every nonempty open set $O$. Using the fact that $\oi^\kappa$ is dense in $X_\kappa$, one shows easily that $\Pi$ is topologically transitive on $X_\kappa$ iff it is topologically transitive on the $\kappa$-cube. By standard arguments~\cite[Theorem~5.9]{Walters82}, this amounts to the existence of a point in $\oi^\kappa$ (or a $G_\delta$ dense set of such points) whose $\Pi$-orbit is dense.

If $\Pi$ is topologically transitive and the set of points whose $\Pi$-orbit is finite is dense, then we say that $\Pi$ acts \newword{chaotically}. By Lemma~\ref{ref16}(iii) and Corollary~\ref{ref17}, $\Sigma_n$ acts chaotically both on $X_n$ and on $\oi^n$. It is well known that a chaotic action on a space such as $\oi^n$ implies sensitive dependence on initial conditions, hence chaotic behaviour in the sense of Devaney~\cite{Banksetal92}.

One constructs easily chaotic elements of $\Sigma_n$. Indeed, the standard tent map on $\oi$ is a McNaughton function, expressible by the formula $s(x_0)=(x_0\land\neg x_0)\oplus(x_0\land\neg x_0)$.
The substitution $\sigma:x_i\mapsto s(x_i)$, for $i<n$, induces therefore on $\oi^n$ the direct product of~$n$ tent maps, which is mixing w.r.t.~Lebesgue measure. Hence $\sigma$ acts in a  topologically transitive and chaotic way. It is not so easy to construct elements of $\Xi_n$ which are chaotic; we will obtain such mappings for even $n$ in Corollary~\ref{ref20}.

\begin{lemma}
Let $\sigma\in\Xi_n$, let $S$ be the induced McNaughton homeomorphism of $\oi^n$, and let $W$ be a complex over the $n$-cube such that $S$ has the form~$(**)$ on each $n$-dimensional cell $C_j$ of $W$. Then all the matrices $A_j$ have the same determinant, which is either $+1$ or $-1$.
\end{lemma}
\begin{proof}
Since the inverse of $S$ is expressible as in $(**)$ via matrices and vectors having integer entries, it is clear that all matrices $A_j$ are invertible and their inverses have integer entries. Therefore all $A_j$'s have determinant $\pm1$. Suppose by contradiction that the $n$-dimensional cells of $W$ are $\vect C1k$ and that there is some $1\le r<k$ such that $\det(A_j)=+1$ for $1\le j\le r$, and $\det(A_j)=-1$ 
for $r<j\le k$. If $1\le j'\le r$ and $r<j''\le k$, then $C_{j'}$ and 
$C_{j''}$ cannot intersect in an $(n-1)$-dimensional face, because 
this would contradict the injectivity of $S$. Let $D$ be the 
topological interior of the $n$-cube,
$E=D\cap(C_1\cup\cdots\cup C_r)$,
$F=D\cap(C_{r+1}\cup\cdots\cup C_k)$. 
Then $E$ and $F$ are nonempty, and closed in the relative topology of 
$D$. By~\cite[Theorem~1.17]{kelley}, $(E\setminus F)\cup(F\setminus 
E)=D\setminus(E\cap F)$
is not connected, and this contradicts~\cite[Theorem~3.61]{hockingyou}, 
since $E\cap F$ has topological 
dimension $\le n-2$.
\end{proof}

\begin{corollary}
For\label{ref18} every $n$ and every $\sigma\in\Xi_n$, the homeomorphism $S$ preserves the Lebesgue measure on $\oi^n$.
\end{corollary}

\begin{corollary}
The only invertible substitutions on $\free 1\MV$ are the identity and the flip $x_0\mapsto \neg x_0$.
\end{corollary}
\begin{proof}
If $S$ is induced by $\sigma\in\Xi_1$, then either $S$ has the form $x_0+k_p$ (with $k_p\in\Zbb$) in every $p$ in which $S$ is differentiable, or the form $-x_0+k_p$. Since the range of $S$ is $\oi$, it must be $k_p=0$ in the first case, or $k_p=1$ in the second.
\end{proof}

The following is the main result of~\cite{pantibernoulli}.

\begin{theorem}
There\label{ref19} is an explicitly constructible family $\{\sigma_{lm}:1\le l,m\in\Nbb\}$ of elements of $\Xi_2$ such that, for every $\sigma_{lm}$ in the family, the induced McNaughton homeomorphism $S_{lm}$ has the following properties:
\begin{itemize}
\item[(i)] $S_{lm}$ fixes pointwise the boundary of $\oi^2$;
\item[(ii)] $S_{lm}$ is ergodic with respect to the Lebesgue measure of the unit square;
\item[(iii)] $S_{lm}$ is non-uniformly hyperbolic and Bernoulli.
\end{itemize}
\end{theorem}

We recall that a measure-preserving bijection is \newword{Bernoulli} if it is measure-theoretically isomorphic to a Bernoulli $2$-sided full shift.

\begin{corollary}
For~\label{ref20} every even $n$ there exist Bernoulli McNaughton homeomorphisms of $\oi^n$. These mappings are mixing w.r.t.~Lebesgue measure,
topologically transitive, and chaotic in the sense of Devaney.
\end{corollary}
\begin{proof}
The Bernoulli property implies mixing, and is preserved under direct products~\cite[Ch.~10 \S1]{CornfeldFomSi82}. Since nonvoid open subsets of $\oi^n$ have positive Lebesgue measure, mixing homeomorphisms are topologically transitive, and hence chaotic by Corollary~\ref{ref17}.
\end{proof}

It would be very interesting to construct a McNaughton homeomorphism of $\oi^3$ having the Bernoulli property: this would allow to extend Corollary~\ref{ref20} to all $n\ge 2$. Up to now, we can only prove the following result about the action of $\Xi_n$ as a group.

\begin{theorem}
For\label{ref23} every $n\ge2$, $\Xi_n$ acts chaotically on $\oi^n$.
\end{theorem}
\begin{proof}
As discussed above, we have a stronger result for even $n$, so we assume $n$ odd. Let $S$ and $T$ be topologically transitive McNaughton homeomorphisms of $\oi^{n-1}$ and $\oi^2$, respectively.
Consider the following direct products:
\begin{align*}
Q&=S\times(\text{identity map on the last coordinate});\\
R&=(\text{identity map on the first $n-2$ coordinates})
\times T.
\end{align*}
Both $Q$ and $R$ are McNaughton homeomorphisms of $\oi^n$, induced by elements of~$\Xi_n$. For every $i<n$, let $0\le a_i<a'_i\le1$
and $0\le b_i<b'_i\le1$: we want to show that an appropriate composition of $Q$ and $R$ maps some point of the open box $A=\prod_{i<n}(a_i,a'_i)$ in the open box
$B=\prod_{i<n}(b_i,b'_i)$.
Since $S$ is topologically transitive, there exists $h\ge0$ such that the open set
$$
U=Q^h[A]\cap\bigl((b_0,b'_0)\times\cdots\times(b_{n-2},b'_{n-2})
\times(a_{n-1},a'_{n-1})\bigr)
$$
is nonempty. Let $0\le c_i<c'_i\le1$ be such that the box
$C=\prod_{i<n}(c_i,c'_i)$ is contained in $U$. The open box
$D=\bigl(\prod_{i<n-1}(c_i,c'_i)\bigr)\times(b_{n_1},b'_{n-1})$
is contained in $B$ and, since $T$ is topologically transitive,
there exists $k\ge0$ such that $R^k[C]\cap D\not=\emptyset$. Therefore
$(R^k\circ Q^h)[A]\cap B\not=\emptyset$.
\end{proof}

For every Borel probability measure $\mu$ on $\oi^n$, let $f_\mu(r)$ denote the integral of the formula $r(\vect x0{{n-1}})$ (viewed as a function $r:\oi^n\to\oi$) with respect to $\mu$. The number
$f_\mu(r)$ may be thought of as the ``average truth-value'' of $r$ w.r.t.~$\mu$. It is natural to restrict attention to measures which are \newword{faithful} ($r\not=0$ implies $f_\mu(r)\not=0$) and \newword{automorphism-invariant} ($f_\mu(r)=f_\mu(\sigma(r))$, for every $\sigma\in\Xi_n$). This latter property is particularly relevant: it says that the average truth-value of a formula should be intrinsic to the formula, and not depending on the particular embedding of it in $\free n\MV$. Lebesgue measure $\lambda$ is faithful (clearly) and automorphism-invariant (by Corollary~\ref{ref18}). 
Let $\mu$ be another probability measure on $\oi^n$, absolutely continuous with respect to $\lambda$. We may assume that $n$ is even, possibly introducing a dummy variable. By Corollary~\ref{ref20}, there exists an automorphism $\sigma$ of
$\free n\MV$ that induces a mixing homeomorphism $S$ on $\oi^n$.
Therefore
the push forward $S_*^k\mu$ of $\mu$ by $S^k$ converges to $\lambda$ in the weak${}^*$ topology~\cite[\S4.9 and Theorem~6.12(ii)]{Walters82}. In particular, for $r$ as above we get
$$
\lim_{k\to\infty}f_\mu(\sigma^k(r))=
\lim_{k\to\infty}\int r\circ S^k\,d\mu=
\lim_{k\to\infty}\int r\,d(S_*^k)\mu=
\int r\,d\lambda=f_\lambda(r).
$$
Hence the existence of mixing McNaughton homeomorphisms gives a distinguished status to $\lambda$. It appears plausible that the only ergodic $\Xi_n$-invariant measures on $\oi^n$ are $\lambda$ and the measures supported on finite orbits. We leave this as an open problem: since measures supported on finite orbits are not faithful, a positive answer would imply that the only reasonable averaging measure on truth-values in \Luk\ logic is Lebesgue measure.

\newcommand{\noopsort}[1]{}

\end{document}